\documentclass[11pt]{article}
\usepackage{anysize}
\marginsize{2.0cm}{2.0cm}{2.0 cm}{2.0 cm}
\usepackage{setspace}

\usepackage{latexsym,amsfonts,amsmath}
\usepackage{graphicx}
\usepackage{color}
\usepackage{epstopdf}
\usepackage[titletoc]{appendix}

\def\liminf{\mathop{\underline{\lim}}}

\newtheorem{condition}{Condition}[section]{\bfseries}{\itshape}

{\bfseries}{\itshape}

\newtheorem{theorem}{Theorem}[section]{\bfseries}{\itshape}

\newtheorem{corollary}{Corollary}[section]{\bfseries}{\itshape}

\newtheorem{proposition}{Proposition}[section]{\bfseries}{\itshape}

{\bfseries}{\itshape}

\newtheorem{lemma}{Lemma}[section]{\bfseries}{\itshape}

\newtheorem{remark}{Remark}[section]{\bfseries}{\itshape}

{\bfseries}{\itshape}

{\bfseries}{\itshape}

\begin{document}
\title{On risk-sensitive piecewise deterministic Markov decision processes}
\author{Xin Guo \thanks{Department of Mathematical Sciences, University of
Liverpool, Liverpool, L69 7ZL, U.K.. E-mail:  X.Guo21@liv.ac.uk.}~ and Yi
Zhang \thanks{Corresponding author. Department of Mathematical Sciences, University of
Liverpool, Liverpool, L69 7ZL, U.K.. E-mail: yi.zhang@liv.ac.uk.}}
\date{}
\maketitle

\par\noindent{\bf Abstract:}  We consider a piecewise deterministic Markov decision process, where the expected exponential utility of total (nonnegative) cost is to be minimized. The cost rate, transition rate and post-jump distributions are under control. The state space is Borel, and the transition and cost rates are locally integrable along the drift. Under natural conditions, we establish the optimality equation, justify the value iteration algorithm, and show the existence of a deterministic stationary optimal policy. Applied to special cases, the obtained results already significantly improve some existing results in the literature on finite horizon and infinite horizon discounted risk-sensitive continuous-time Markov decision processes.
\bigskip

\par\noindent {\bf Keywords:} Continuous-time Markov decision
processes. Piecewise deterministic Markov decision processes. Exponential utility. Dynamic programming.
\bigskip

\par\noindent
{\bf AMS 2000 subject classification:} Primary 90C40,  Secondary
60J75

\section{Introduction}
 Since the pioneering work \cite{Howard:1972}, risk-sensitive discrete-time Markov decision processes (DTMDPs) have been studied intensively. Having restricted our attention to total undiscounted or discounted problems, let us mention e.g., \cite{BauerleRieder:2014,Cavazos:2000,Chung:1987,Di Masi:1999,Fainberg:1982,Jaquette:1976,Jaskiewicz:2008}, most of which deal with the exponential utility, as well as in the present paper. As an application, an open problem in insurance was recently solved in \cite{BauerleAnna:2015} in the framework of risk-sensitive DTMDP. There are notable differences between risk-sensitive and risk-neutral DTMDPs. For instance, in a finite model, i.e., when the state and action spaces are both finite, there is always a deterministic stationary optimal policy in a discounted  risk-neutral DTMDP, but not always in a discounted risk-sensitive DTMDP, see \cite{Jaquette:1976}.

One of the first works on risk-sensitive continuous-time Markov decision processes (CTMDPs) is \cite{Piunovski:1985}, where only verification theorems were presented. Recently, there have been reviving interests in this topic; see e.g., \cite{Coraluppi:1997,Ghosh:2014,Kumar:2013,Wei:2016,WeiC:2016,Zhang:2017}. A finite horizon total undiscounted risk-sensitive CTMDP was considered in \cite{Ghosh:2014,Piunovski:1985,Wei:2016}, whose arguments were summarized as follows. Firstly, the optimality equation is shown to admit a solution out of a small enough class. Secondly, by using the Feynman-Kac formula, this solution is shown to be the value function, and any Markov policy providing the minimizer in the optimality equation is optimal. The proofs of \cite{Ghosh:2014,Wei:2016} reveal that the the main technicalities lie in the first step, for which, the state space was assumed to be denumerable. This assumption is important for the diagonalization argument used in \cite{Wei:2016}, which is an extension of \cite{Ghosh:2014} from bounded transition rate to possibly unbounded transition rate, whose growth is bounded by a Lyapunov function. The latter requirement and the boundedness of the cost rate then validate the Feynman-Kac formula applied in the second step. The author of \cite{Wei:2016} mentioned that it was unclear how to extend his argument to an unbounded cost rate, see Section 7 therein. Following a similar argument as described above, a discounted risk-sensitive CTMDP was also considered in \cite{Ghosh:2014}, although now the first step becomes, to quote the authors' words (see p.658 therein), ``surprisingly far more involved'', for which the state space was further assumed to be finite, see Remark 3.6 therein. It is a corollary of the present paper that we significantly weaken the restrictive conditions in \cite{Ghosh:2014,Wei:2016}, see Section \ref{GGZyPDsEC3} below.

The present paper is concerned with a risk-sensitive piecewise deterministic Markov decision process (PDMDP), where the expected exponential utility of the total cost is to be minimized. The state space is a general Borel space, the transition and the nonnegative cost rates only need be locally integrable along the drift. A PDMDP is an extension of a CTMDP: now between two consecutive jumps, the process evolves according to deterministic Markov process. For simplicity and to keep the conditions as weak as possible, we do not consider the control on the drift. In spite that there has been a vast literature on PDMDPs; see the well known monographs \cite{Costa:2013,Davis:1993} and the references therein, to the best of our knowledge, risk-sensitive PDMDPs have not been systematically studied before.

Our main contributions are the following. We establish the optimality equation satisfied by the value function, justify the value iteration algorithm and show the existence of a deterministic stationary optimal policy. As an application and corollary, finite horizon and infinite horizon discounted risk-sensitive CTMDPs are reformulated as total undiscounted risk-sensitive PDMDPs, and are thus treated in a unified way and under much weaker conditions than in \cite{Ghosh:2014,Wei:2016}. This is possible because we follow a different argument. Namely, we directly show that the value function satisfies the optimality equation, by reducing the total undiscounted risk-sensitive PDMDP to a risk-sensitive DTMDP. This method, without referring to the Feynman-Kac formula, was originally developed by Yushkevich \cite{Yushkevich:1980} for risk-neutral CTMDPs. Later, it was employed in \cite{BauerleRieder:2009,BauerleRieder:2011,Costa:2013,Davis:1993,Forwick:2004,Schal:1998} for studies of risk-neutral PDMDPs, and in \cite{Zhang:2017} for risk-sensitive CTMDPs. In \cite{Coraluppi:1997}, restricted to stationary policies, the discounted risk-sensitive CTMDP with bounded transition rates was reduced to a DTMDP problem, using the uniformization technique. The induced DTMDP is less standard (with a random cost), and was not further investigated there.

The rest of the paper is organized as follows. In Section \ref{GGZYPDMDPSec2} we describe the concerned optimal control problem. In Section \ref{GGZyPDsEC3} we present the main results, the proofs of which are postponed to Section \ref{GGZYPDMDSec5}. We finish the paper with a conclusion in Section \ref{GGZPDMDPSecConclusion}. Some relevant facts were collected in the appendix for ease of reference.

\section{Model description and problem statement}\label{GGZYPDMDPSec2}
\par\noindent\textbf{Notations and conventions.} In what follows, ${\cal{B}}(X)$ is
the Borel $\sigma$-algebra of the topological space $X,$ $I$ stands for the indicator function, and $\delta_{\{x\}}(\cdot)$
is the Dirac measure concentrated on the singleton $\{x\},$ assumed to be measurable. A measure is $\sigma$-additive and $[0,\infty]$-valued. Below, unless stated otherwise, the term of
measurability is always understood in the Borel sense. Throughout
this paper, we adopt the conventions of
\begin{eqnarray}\label{GGZyExponential56}
\frac{0}{0}:=0,~0\cdot\infty:=0,~\frac{1}{0}:=+\infty,~\infty-\infty:=\infty.
\end{eqnarray}
If a mapping $f$ defined on $X$, and $\{X_i\}$ is a partition of $X$, then when $f$ is piecewise defined as $f(x)=g_i(x)$ for all $x\in X_i$, the notation $f(x)=\sum_i I\{x\in X_i\}g_i(x)$ is used, even if $f$ is not real-valued.

Let $S$ be a nonempty Borel state space, $A$ be
a nonempty Borel action space, and $q$
stand for a signed kernel $q(dy|x,a)$ on ${\cal{B}}(S)$ given
$(x,a)\in S\times A$ such that
\begin{eqnarray}\label{GGZYPDMDP20}
\tilde{q}(\Gamma_S|x,a):=q(\Gamma_S\setminus\{x\}|x,a)\ge 0
\end{eqnarray}
for all $\Gamma_S\in{\cal{B}}(S).$ Throughout this article we assume
that $q(\cdot|x,a)$ is conservative and stable, i.e.,
\begin{eqnarray}\label{GGZZPDMDP19}
q(S|x,a)=0,~\bar{q}_x=\sup_{a\in A}q_x(a)<\infty,
\end{eqnarray}
where $q_x(a):=-q(\{x\}|x,a).$ The signed kernel $q$ is often called the
transition rate. Between two consecutive jumps, the state of the process evolves according to a measurable mapping $\phi$ from $S\times [0,\infty)$ to $S$, see (\ref{GGGAMO01}) below. It is assumed that for each $x\in S$
\begin{eqnarray}\label{GGZyExponential01}
\phi(x,t+s)=\phi(\phi(x,t),s),~\forall~s,t\ge 0;~\phi(x,0)=x,
\end{eqnarray}
and $t\rightarrow \phi(x,t)$ is continuous.

Finally let the cost rate $c$ be a $[0,\infty)$-valued measurable function on $S\times A$.
For simplicity, we do not consider the case of different admissible action spaces at
different states.

\begin{condition}\label{GZyExponentialCondition02}
\begin{itemize}
\item[(a)] For each bounded measurable function $f$ on $S$ and each $x\in S$, $\int_S f(y)\tilde{q}(dy|x,a)$ is continuous in $a\in A.$
\item[(b)] For each $x\in S,$ the (nonnegative) function $c(x,a)$ is lower semicontinuous in $a\in A.$
\item[(c)] The action space $A$ is a compact Borel space.
\end{itemize}
\end{condition}

\begin{condition}\label{GGZyExponentialConditionExtra}
For each $x\in S$, $\int_{0}^t\overline{q}_{\phi(x,s)}ds<\infty$, and $\int_{0}^t \sup_{a\in A} c(\phi(x,s),a)ds<\infty$, for each $t\in[0,\infty).$
\end{condition}
The integrals in the above condition are well defined: the integrands are universally measurable in $s\in[0,\infty)$; see Chapter 7 of \cite{Bertsekas:1978}.

Let us take the sample space $\Omega$ by adjoining to the
countable product space $S\times((0,\infty)\times S)^\infty$ the
sequences of the form
$(x_0,\theta_1,\dots,\theta_n,x_n,\infty,x_\infty,\infty,x_\infty,\dots),$
where $x_0,x_1,\dots,x_n$ belong to $S$,
$\theta_1,\dots,\theta_n$ belong to $(0,\infty),$ and
$x_{\infty}\notin S$ is the isolated point. We equip $\Omega$ with
its Borel $\sigma$-algebra $\cal F$.

Let $t_0(\omega):=0=:\theta_0,$ and for each $n\geq 0$, and each
element $\omega:=(x_0,\theta_1,x_1,\theta_2,\dots)\in \Omega$, let
\begin{eqnarray*}
t_n(\omega)&:=&t_{n-1}(\omega)+\theta_n,
\end{eqnarray*}
and
\begin{eqnarray*}
t_\infty(\omega):=\lim_{n\rightarrow\infty}t_n(\omega).
\end{eqnarray*}
Obviously, $(t_n(\omega))$ are measurable mappings on $(\Omega,{\cal
F})$. In what follows, we often omit the argument $\omega\in
\Omega$ from the presentation for simplicity. Also, we regard
$x_n$ and $\theta_{n+1}$ as the coordinate variables, and note
that the pairs $\{t_n,x_n\}$ form a marked point process with the
internal history $\{{\cal F}_t\}_{t\ge 0},$ i.e., the filtration
generated by $\{t_n,x_n\}$; see Chapter 4 of \cite{Kitaev:1995}
for greater details. The marked point process $\{t_n,x_n\}$
defines the stochastic process $\{\xi_t,t\ge 0\}$ on $(\Omega,{\cal F})$ of interest
by
\begin{eqnarray}\label{GGGAMO01}
\xi_t=\sum_{n\ge 0}I\{t_n\le t<t_{n+1}\}\phi(x_n,t-t_n)+I\{t_\infty\le
t\}x_\infty,~t\ge 0,
\end{eqnarray}
where we accept $0\cdot x:=0$ and $1\cdot x:=x$ for each $x\in S_\infty,$ and below we denote
$S_{\infty}:=S\bigcup\{x_\infty\}$.

A (history-dependent) policy $\pi$ is given by a
sequence $(\pi_n)$ such that, for each $n=0,1,2,\dots,$
$\pi_n(da|x_0,\theta_1,\dots,x_{n},s)$ is a stochastic kernel on
$A$, and for each $\omega=(x_0,\theta_1,x_1,\theta_2,\dots)\in
\Omega$, $t> 0,$
\begin{eqnarray}\label{GGZZPDMDPStarNumber}
\pi(da|\omega,t)&=&I\{t\ge t_\infty\}\delta_{a_\infty}(da)+
\sum_{n=0}^\infty I\{t_n< t\le
t_{n+1}\}\pi_{n}(da|x_0,\theta_1,\dots,\theta_n,x_n, t-t_n),
\end{eqnarray}
where $a_\infty\notin A$ is some isolated point. A policy
$\pi$  is called Markov if, with slight abuse of
notations,  $
\pi(da|\omega,s)=\pi^M(da|\xi_{s-},s)$ for some stochastic kernel $\pi^M$. A
Markov policy is further called deterministic if the stochastic
kernels $\pi^M(da|x,s)=\delta_{\{f^M(x,s)\}}(da)$ for some measurable mapping $f^M$ from $S\times(0,\infty)$ to $A$.  A policy is called deterministic stationary if for each $n=0,1,\dots,$ $\pi_{n}(da|x_0,\theta_1,\dots,\theta_n,x_n, t-t_n)=\delta_{\{f(\phi(x_n,t-t_n))\}}(da)$
for some measurable mapping $f$ from $S$ to $A$. We shall identify such a deterministic stationary policy by the underlying measurable mapping $f$.

The class of all policies is denoted by $\Pi.$
Under a fixed policy $\pi=(\pi_n)$, for each initial distribution $\gamma$ on $(S,{\cal B}(S)),$ by using the Ionescu-Tulcea theorem, one can build a probability measure $P_\gamma^\pi$ on $(\Omega,{\cal F})$ such that $P_\gamma^\pi(x_0\in \Gamma)=\gamma(\Gamma)$ for each $\Gamma\in {\cal B}(S)$, and the conditional distribution of $(\theta_{n+1},x_{n+1})$ with the condition on $x_0,\theta_1,x_1,\dots,\theta_{n},x_n$ is given on $\{\omega:x_n(\omega)\in S\}$ by
\begin{eqnarray}\label{GGZyExponential21}
&&P_\gamma^\pi(\theta_{n+1}\in \Gamma_1,~x_{n+1}\in \Gamma_2|x_0,\theta_1,x_1,\dots,\theta_{n},x_n)\nonumber\\
&=&\int_{\Gamma_1}e^{-\int_0^t \int_A q_{\phi(x_n,s)}(a)\pi_n(da|x_0,\theta_1,\dots,\theta_n,x_n,s)ds}\int_{A}\tilde{q}(\Gamma_2|\phi(x_n,t),a)\pi_n(da|x_0,\theta_1,\dots,\theta_n,x_n,t)dt,\nonumber\\
&&~\forall~\Gamma_1\in{\cal B}((0,\infty)),~\Gamma_2\in{\cal B}(S);\nonumber\\
&&P_\gamma^\pi(\theta_{n+1}=\infty,~x_{n+1}=x_\infty|x_0,\theta_1,x_1,\dots,\theta_{n},x_n)=e^{-\int_0^\infty  \int_A q_{\phi(x_n,s)}(a)\pi_n(da|x_0,\theta_1,\dots,\theta_n,x_n,s)ds},
\end{eqnarray}
and given on $\{\omega:x_n(\omega)=x_\infty\}$ by
\begin{eqnarray*}
P_\gamma^\pi(\theta_{n+1}=\infty,~x_{n+1}=x_\infty|x_0,\theta_1,x_1,\dots,\theta_{n},x_n)=1.
\end{eqnarray*}
Below, when $\gamma$ is a Dirac measure concentrated at $x\in S,$
we use the denotation ${}{P}_x^\pi.$ Expectations with respect to
${}{P}_\gamma^\pi$ and ${}{P}_x^\pi$ are denoted as
${}{E}_{\gamma}^\pi$ and ${}{E}_{x}^\pi,$ respectively. Roughly speaking, the uncontrolled version of the process evolves as follows: given the current state, the process evolves deterministically according to the mapping $\phi$, up to the next jump, taking place after a random time whose distribution is (nonstationary) exponential, and the dynamics continue in the similar manner. A detailed book treatment with many examples of this and more general type of processes, allowing deterministic jumps, can be found in \cite{Davis:1993}.

For each $x\in S$, and policy $\pi=(\pi_n)$,
\begin{eqnarray*}
E_x^\pi\left[e^{\int_0^\infty \int_A
c(\xi_t,a)\pi(da|\omega,t)dt}\right]=E_x^\pi\left[e^{\sum_{n=0}^\infty \int_0^{\theta_{n+1}} \int_A
c(\phi(x_n,s),a)\pi_n(da|x_0,\theta_1,\dots,x_n,s)ds}\right]=:V(x,\pi)
\end{eqnarray*}
defines the concerned performance measure of the policy $\pi\in \Pi$ given the initial state $x\in S.$
Here and below, we put $c(x_\infty,a):=0$ for each $a\in A,$ and $\phi(x_\infty,t)=x_\infty$ for each $t\in [0,\infty).$ We are interested in the following optimal control problem for each $x\in S:$
\begin{eqnarray}\label{GZPDMDP17}
&\mbox{Minimize over $\pi\in\Pi$:~}&V(x,\pi).
\end{eqnarray}
A policy $\pi^\ast$ is called optimal if
$
V(x,\pi^\ast)=\inf_{\pi\in \Pi}V(x,\pi)=:V^\ast(x)$ for each $x\in S$.

The objective of this paper is to show, under the imposed conditions, the existence of a deterministic stationary optimal policy, and to establish the corresponding optimality equation satisfied by the value function $V^\ast$, together with its value iteration.
Evidently, $V^\ast(x)\ge 1$ for each $x\in S.$  Under the next condition, it will be seen that for each $x\in S,$ $V^\ast(\phi(x,s))$ is absolutely continuous in $s.$

\begin{condition}\label{GGZyExponentialFinitenessCon}
For each $x\in S,$ $V^\ast(x)<\infty$.
\end{condition}
The above condition is mainly assumed for notational convenience. In fact, the main optimality results (such as the existence of a deterministic stationary optimal policy) obtained in this paper can be established without assuming Condition \ref{GGZyExponentialFinitenessCon}, at the cost of some additional notations. In a nutshell, one has to consider the sets
$
\hat{S}:=\{x\in S:~V^\ast(x)<\infty\}
$
and $S\setminus \hat{S}$ separately, and note that if $x\in \hat{S}$, then $\phi(x,t)\in \hat{S}$ for each $t\in [0,\infty).$ The reasoning presented under Condition \ref{GGZyExponentialFinitenessCon} can be followed in an obvious manner. We formulate the corresponding optimality results in Remarks \ref{GGGAMORem01} and \ref{GGGAMORem02} below.

\section{Main statements}\label{GGZyPDsEC3}
We first present the main optimality results concerning problem (\ref{GZPDMDP17}) for the PDMDP model. Their proofs are postponed to the next section.
 \begin{theorem}\label{GGZyExponentialTheorem}
 Suppose Conditions \ref{GZyExponentialCondition02}, \ref{GGZyExponentialConditionExtra} and \ref{GGZyExponentialFinitenessCon} are satisfied. Then the following assertions hold.
 \begin{itemize}
 \item[(a)] The value function $V^\ast$ for problem (\ref{GZPDMDP17}) is the minimal $[1,\infty)$-valued solution to the following optimality equation:
  \begin{eqnarray*}
&& -(V(\phi(x,t))-V(x))\\
&=&\int_0^t \inf_{a\in A}\left\{ \int_S V(y)\tilde{q}(dy|\phi(x,\tau),a)- (q_{\phi(x,\tau)}(a)-c(\phi(x,\tau),a) )V(\phi(x,\tau))\right\}d\tau,\\
&&~t\in[0,\infty),x\in S.
 \end{eqnarray*}
In particular, $V^\ast(\phi(x,t))$ is absolutely continuous in $t$ for each $x\in S.$
 \item[(b)]  There exists a deterministic stationary optimal policy $f$, which can be taken as any measurable mapping from $S$ to $A$ such that
 \begin{eqnarray*}
 &&\inf_{a\in A}\left\{ \int_S V^\ast(y)\tilde{q}(dy|x,a)- (q_{x}(a)-c(x,a))V^\ast(x))\right\}\\
 &=&\int_S V^\ast(y)\tilde{q}(dy|x,f(x))- (q_{x}(f(x))-c(x,f(x)))V^\ast(x)),~\forall~x\in S.
 \end{eqnarray*}
\end{itemize}
 \end{theorem}

\begin{remark}\label{GGGAMORem01}
By inspecting its proof, one can see the following version of Theorem \ref{GGZyExponentialTheorem} holds without assuming Condition \ref{GGZyExponentialFinitenessCon}. Suppose Conditions \ref{GZyExponentialCondition02} and \ref{GGZyExponentialConditionExtra} are satisfied. Then the following assertions hold.
 \begin{itemize}
 \item[(a)] The value function $V^\ast$ for problem (\ref{GZPDMDP17}) is the minimal $[1,\infty]$-valued solution to the following optimality equation:
  \begin{eqnarray*}
&& -(V(\phi(x,t))-V(x))\\
&=&\int_0^t \inf_{a\in A}\left\{ \int_S V(y)\tilde{q}(dy|\phi(x,\tau),a)- (q_{\phi(x,\tau)}(a)-c(\phi(x,\tau),a) )V(\phi(x,\tau))\right\}d\tau,\\
&&~t\in[0,\infty),x\in \hat{S};\\
&&V(x)<\infty,~x\in \hat{S};~V(x)=\infty,~x\in S\setminus \hat{S}.
 \end{eqnarray*}
In particular, $V^\ast(\phi(x,t))$ is absolutely continuous in $t$ for each $x\in \hat{S}.$
 \item[(b)]  There exists a deterministic stationary optimal policy $f$, which can be taken as any measurable mapping from $S$ to $A$ such that
 \begin{eqnarray*}
 &&\inf_{a\in A}\left\{ \int_S V^\ast(y)\tilde{q}(dy|x,a)- (q_{x}(a)-c(x,a))V^\ast(x))\right\}\\
 &=&\int_S V^\ast(y)\tilde{q}(dy|x,f(x))- (q_{x}(f(x))-c(x,f(x)))V^\ast(x)),~\forall~x\in \hat{S}.
 \end{eqnarray*}
\end{itemize}
\end{remark}

Next, we present the value iteration algorithm for the value function $V^\ast$.
 \begin{theorem}\label{GGGZyExponentialTheorem2}
 Suppose Conditions \ref{GZyExponentialCondition02}, \ref{GGZyExponentialConditionExtra} and \ref{GGZyExponentialFinitenessCon} are satisfied. Let $V^{(0)}(x):=1$ for each $x\in S$. For each $n\ge 0,$ let $V^{(n+1)}$ be the minimal $[1,\infty)$-valued measurable solution to
 \begin{eqnarray}\label{GZYPDMDPEqn11}
 && -(V^{(n+1)}(\phi(x,t))-V^{(n+1)}(x))\nonumber\\
&=&\int_0^t \inf_{a\in A}\left\{ \int_S V^{(n)}(y)\tilde{q}(dy|\phi(x,\tau),a)- (q_{\phi(x,\tau)}(a)-c(\phi(x,\tau),a) )V^{(n+1)}(\phi(x,\tau))\right\}d\tau,\nonumber\\
&&~t\in[0,\infty),x\in S,
 \end{eqnarray}
such that $V^{(n+1)}(\phi(x,t))$ is absolutely continuous in $t$ for each $x\in S.$ (For each $n\ge 0,$ such a solution always exists.) Furthermore, $\{V^{(n)}\}$ is a monontone nondecreasing sequence of measurable functions on $S$ such that for each $x\in S,$ $V^{(n)}(x)\uparrow V^\ast(x)$ as $n\uparrow \infty.$
 \end{theorem}

\begin{remark}\label{GGGAMORem02}
Similar to Remark \ref{GGGAMORem01}, we have the following version of Theorem \ref{GGGZyExponentialTheorem2} without assuming Condition \ref{GGZyExponentialFinitenessCon}. Suppose Conditions \ref{GZyExponentialCondition02}, \ref{GGZyExponentialConditionExtra} are satisfied. Let $V^{(0)}(x):=1$ for each $x\in \hat{S}$ and $V^{(0)}(x)=\infty$ if $x\in S\setminus \hat{S}$. For each $n\ge 0,$ let $V^{(n+1)}$ be the minimal $[1,\infty]$-valued measurable solution to
 \begin{eqnarray*}
 && -(V^{(n+1)}(\phi(x,t))-V^{(n+1)}(x))\nonumber\\
&=&\int_0^t \inf_{a\in A}\left\{ \int_S V^{(n)}(y)\tilde{q}(dy|\phi(x,\tau),a)- (q_{\phi(x,\tau)}(a)-c(\phi(x,\tau),a) )V^{(n+1)}(\phi(x,\tau))\right\}d\tau,\nonumber\\
&&~t\in[0,\infty),x\in \hat{S},\\
&&V^{(n+1)}(x)<\infty,~x\in \hat{S},~V^{(n+1)}(x)=\infty,~x\in S\setminus \hat{S}.
 \end{eqnarray*}
Here $V^{(n+1)}(\phi(x,t))$ is absolutely continuous in $t$ for each $x\in \hat{S}.$ (For each $n\ge 0,$ such a solution always exists.) Furthermore, $\{V^{(n)}\}$ is a monontone nondecreasing sequence of measurable functions on $S$ such that for each $x\in S,$ $V^{(n)}(x)\uparrow V^\ast(x)$ as $n\uparrow \infty.$
\end{remark}

We can apply our theorems to a special case of a CTMDP.  That is, $\phi(x,t)\equiv x$ for each $x\in S.$ The following $\alpha$-discounted risk-sensitive CTMDP problem was considered in \cite{Ghosh:2014}:
\begin{eqnarray}\label{GGZYPDMDP16}
\mbox{Minimize over $\pi\in \Pi$: } E_x^\pi\left[e^{\int_0^\infty e^{-\alpha t}\int_A
c(\xi_t,a)\pi(da|\omega,t)dt}\right],~x\in S.
\end{eqnarray}
Here $\alpha>0$ is a fixed constant.
In fact, the authors of \cite{Ghosh:2014} were restricted to Markov policies, bounded transition and cost rates, i.e., $\sup_{x\in S}\overline{q}_x<\infty$, and $\sup_{x\in S,a\in A}c(x,a)<\infty$, and a finite state space $S$. These restrictions, e.g., the finiteness of $S$, were needed for their investigations, see e.g., Remark 3.6 in \cite{Ghosh:2014}. Under the compactness-continuity condition (Condition \ref{GZyExponentialCondition02}), it was shown in \cite{Ghosh:2014} that there exists an optimal Markov policy for the discounted risk-sensitive CTMDP, and established the optimality equation. By using the theorems presented earlier in this section, we can obtain these optimality results for problem (\ref{GGZYPDMDP16}) in a much more general setup: the state space $S$ is Borel, there is no boundedness requirement on the transition rate with respect to the state $x\in S$, and the optimality is over the class of history-dependent policies. Furthermore, we let the CTMDP model be nonhomogeneous, i.e., the transition rate $q(dy|t,x,a)$ now is a signed kernel on ${\cal B}(S)$ from $(t,x,a)\in[0,\infty)\times S\times A$, satisfying the corresponding version of (\ref{GGZZPDMDP19}); the notations $\tilde{q}$ is kept as before, see (\ref{GGZYPDMDP20}), with the extra argument $t$ in addition to $x$. Similarly, the nonnegative cost rate $c$ is allowed to be a measurable function on $[0,\infty)\times S\times A$.

\begin{corollary}
Consider the $\alpha$-discounted risk-sensitive (nonhomogeneous) CTMDP problem (\ref{GGZYPDMDP16}) with $c(\xi_t,a)$ being replaced by $c(t,\xi_t,a)$. Suppose
\begin{eqnarray*}
\sup_{t\in[0,\infty)}\{\overline{q}_{(t,x)}\}<\infty,~\forall~x\in S,~ \sup_{t\in[0,\infty),x\in S,a\in A}c(t,x,a)<\infty,
\end{eqnarray*}
and the corresponding version of Condition \ref{GZyExponentialCondition02}, where $x$ is replaced by $(t,x)$, is satisfied by the nonhomogeneous CTMDP model. Then the following assertions hold.
\begin{itemize}
\item[(a)] There exists some $[1,\infty)$-valued measurable solution on $[0,\infty)\times S$ to
\begin{eqnarray*}
&&-(V(t,x)-V(0,x))\\
&=&\int_0^t \inf_{a\in A}\left\{\int_S V(u,y)\tilde{q}(dy|u,x,a)+(e^{-\alpha u} c(u,x,a)-q_{(u,x)}(a))V(u,x)  \right\}du,\\
&&~x\in S,~t\in [0,\infty),
\end{eqnarray*}
so that $V(t,x)$ is absolutely continuous in $t$ for each $x\in S.$
\item[(b)] Let $L$ be the minimal $[1,\infty)$-valued measurable solution on $[0,\infty)\times S$ to the above equation. Then
the value function say $L^\ast$ to the $\alpha$-discounted risk-sensitive CTMDP problem (\ref{GGZYPDMDP16}) (with $c(\xi_t,a)$ being replaced by $c(t,\xi_t,a)$) is given by $L^\ast(x)=L(0,x)$ for each $x\in S.$

\item[(c)] There exists an optimal deterministsic Markov policy $f$ for the $\alpha$-discounted risk-sensitive CTMDP problem (\ref{GGZYPDMDP16}) (with $c(\xi_t,a)$ being replaced by $c(t,\xi_t,a)$). One can take $f$ as any measurable mapping from $[0,\infty)\times S$ to $A$ such that
    \begin{eqnarray*}
    &&\inf_{a\in A}\left\{\int_S L(u,y)\tilde{q}(dy|u,x,a)+(e^{-\alpha u} c(u,x,a)-q_{(u,x)}(a))L(u,x)  \right\}\\
    &=&\int_S L(u,y)\tilde{q}(dy|u,x,f(u,x))+(e^{-\alpha u} c(u,x,f(u,x))-q_{(u,x)}(f(u,x)))L(u,x)
    \end{eqnarray*}
    for each $u\in [0,\infty)$ and $x\in S.$
\end{itemize}
 \end{corollary}

\par\noindent\textit{Proof.} We prove this by reformulating the nonhomogeneous version of the $\alpha$-discounted risk-sensitive (nonhomogeneous) CTMDP problem (\ref{GGZYPDMDP16}) in the form of problem
(\ref{GZPDMDP17}) for a PDMDP, which we introduce as follows. We use the notation ``hat'' to distinguish this model from the original (nonhomogeneous) CTMDP model.

\begin{itemize}
\item The state space is $\hat{S}=[0,\infty)\times S.$
\item The action space is the same as in the CTMDP: $\hat{A}=A.$
\item the transition rate $\hat{q}(ds\times dy|(t,x),a)$ is defined by
\begin{eqnarray*}
\hat{q}(ds\times dy|(t,x),a):=\tilde{\hat{q}}(ds\times dy|(t,x),a)-I\{(t,x)\in ds\times dy\}q_{(t,x)}(a),
\end{eqnarray*}
where \begin{eqnarray*}
\tilde{\hat{q}}(ds\times dy|(t,x),a):=I\{t\in ds\}\tilde{q}(dy|t,x,a),
\end{eqnarray*}
for each $(t,x)\in \hat{S}$ and $a\in  \hat{A}.$
\item The drift is given by $\hat{\phi}((t,x),s):=(t+s,x)$ for each $x\in S$ and $t,s\ge 0.$ Clearly it satisfies the corresponding version of (\ref{GGZyExponential01}).
\item The cost rate is given by
\begin{eqnarray*}
\hat{c}((t,x),a):=e^{-\alpha t} c(t,x,a),~\forall~t\in[0,\infty),~x\in S,~a\in A.
\end{eqnarray*}
\end{itemize}
Now the marked point process $\{\hat{t}_n,\hat{x}_n\}$ and controlled process $\hat{\xi}_t$ in this PDMDP model is connected to those in the original (nonhomogeneous) CTMDP model, namely $(t_n,x_n)$ and $\xi_t$, via  $\hat{t}_n=t_n$ and $\hat{x}_n=(t_n,x_n),$ and $\hat{\xi}_t=(t,\xi_t).$  For example, under a fixed strategy $\hat{\pi}$  and initial distribution $\hat{\gamma}$ in this PDMDP model, the version of the first equation in (\ref{GGZyExponential21}) now reads on $\{\omega:x_n(\omega)\in S\}$
\begin{eqnarray*}
&&\hat{P}_{\hat{\gamma}}^{\hat{\pi}}(\hat{\theta}_{n+1}\in \Gamma_1,~\hat{x}_{n+1}\in \Gamma_2\times \Gamma_3|\hat{x}_0,\hat{\theta}_1,\hat{x}_1,\dots,\hat{\theta}_{n},\hat{x}_n)\nonumber\\
&=&\int_{\Gamma_1}e^{-\int_0^t \int_A q_{(t_n+s,x_n)}(a)\hat{\pi}_n(da|\hat{x}_0,\hat{\theta}_1,\dots,\hat{\theta}_n,\hat{x}_n,s)ds}\\
&&\times \int_{A}I\{t+t_n\in\Gamma_2\} \tilde{q} (\Gamma_3|t+t_n,x_n ,a)\hat{\pi}_n(da|\hat{x}_0,\hat{\theta}_1,\dots,\hat{\theta}_n,\hat{x}_n,t)dt,\nonumber\\
&&~\forall~\Gamma_1\in{\cal B}((0,\infty)),~\Gamma_2\in {\cal B}([0,\infty)), ~\Gamma_3\in{\cal B}(S).\nonumber
\end{eqnarray*}

Clearly, Conditions \ref{GZyExponentialCondition02}, \ref{GGZyExponentialConditionExtra} and \ref{GGZyExponentialFinitenessCon} are satisfied by this PDMDP model. It remains to apply Theorem \ref{GGZyExponentialTheorem}.
$\hfill\Box$
\bigskip

The condition in the previous corollary is much weaker than in \cite{Ghosh:2014}, and can be further weakened; one only needs the reformulated PDMDP to satisfy Conditions \ref{GZyExponentialCondition02}, \ref{GGZyExponentialConditionExtra} and \ref{GGZyExponentialFinitenessCon}. Moreover, the boundedness of the cost rate $c$ was assumed in the previous corollary only to ensure Condition \ref{GGZyExponentialFinitenessCon} to be satisfied. It can be relaxed if one formulates the previous corollary using the statements in Remarks \ref{GGGAMORem01} and \ref{GGGAMORem02}.

One can also consider the risk-sensitive nonhomogeneous CTMDP problem on the finite horizon $[0,T]$ with $T>0$ being a fixed constant:
\begin{eqnarray*}
\mbox{Minimize over $\pi\in \Pi$: } E_x^\pi\left[e^{\int_0^T e^{-\alpha t}\int_A
c(t,\xi_t,a)\pi(da|\omega,t)dt +g(\xi_T)}\right],~x\in S,
\end{eqnarray*}
where $g$ is a $[0,\infty)$-valued measurable function; $g(x)$ represents the terminal cost incurred when $\xi_T=x\in S$. Let us put  $g(x_\infty):=0.$  Here $\alpha$ is a fixed nonnegative finite constant. A simpler version of this problem was considered in \cite{Wei:2016} with $\alpha=0$ and a bounded cost rate, where additional restrictions were put on the growth of the transition rate.
We can reformulate this problem into the PDMDP problem (\ref{GZPDMDP17}) just as in the above. The only difference is that now we put $q_{(t,x)}(a)\equiv 0$ for each $x\in S$ and $t\ge T,$ and introduce the following cost rate for each $x\in S$, $t\ge 0$ and $a\in A:$
\begin{eqnarray*}
\hat{c}((t,x),a)=\left\{\begin{array}{ll}
e^{-\alpha t}c(t,x,a), & \mbox{ if }t\le T;  \\
e^{-(t-T)}g(x) & \mbox{ if } t>T. \end{array}\right.
\end{eqnarray*}


%

\section{Proof of the main statements}\label{GGZYPDMDSec5}

For the rest of this paper, it is convenient to introduce the
following notations. Let $\mathbb{P}(A)$ be the space of
probability measures on ${\cal B}(A)$, endowed with the
standard weak topology. For each $\mu\in \mathbb{P}(A)$,
\begin{eqnarray*}
&& q_x(\mu):=\int_A q_x(a)\mu(da),~\tilde{q}(dy|x,\mu):= \int_A \tilde{q}(dy|x,a)\mu(da),~c(x,\mu):=\int_A c(x,a)\mu(da).
\end{eqnarray*}
Let ${\cal R}$ denote the set of (Borel) measurable mappings
$\rho_t(da)$ from $t\in(0,\infty)\rightarrow \mathbb{P}(A).$ Here,
we do not distinguish two measurable mappings in $t\in
(0,\infty),$ which coincide almost everywhere with respect to the
Lebesgue measure. Let us equip ${\cal R}$ with the Young topology, which is the weakest
topology with respect to which the function
$
\rho\in {\cal{R}}\rightarrow \int_0^\infty \int_A
f(t,a)\rho_t(da)dt
$
is continuous for each strongly integrable Carath\'{e}odory function
$f$ on $(0,\infty)\times A$ . Here a real-valued measurable
function $f$ on $(0,\infty)\times A$ is called a strongly
integrable Carath\'{e}odory function if for each fixed
$t\in(0,\infty)$, $f(t,a)$ is continuous in $a\in A,$ and for each
fixed $a\in A,$ $\sup_{a\in A}|f(t,a)|$ is integrable in $t$,
i.e., $\int_0^\infty \sup_{a\in A}|f(t,a)|dt<\infty.$ It is known that if $A$ is a compact Borel space, then so is ${\cal
R}$; see Chapter 4 of \cite{Davis:1993}.

\begin{lemma}\label{GGZyPDMDPLem01}
 Suppose Conditions \ref{GZyExponentialCondition02} and \ref{GGZyExponentialConditionExtra} are satisfied. Then the following assertions hold.
\begin{itemize}
\item[(a)] The value function
$V^\ast$ is the minimal $[1,\infty]$-valued measurable solution to
\begin{eqnarray*}
V^\ast(x)&=& \inf_{\rho\in {\cal R}}\left\{\int_0^\infty e^{-\int_0^\tau (q_{\phi(x,s)}(\rho_s)-c(\phi(x,s),\rho_s))ds} \left(\int_S V^\ast(y)\tilde{q}(dy|\phi(x,\tau),\rho_\tau)\right)d\tau\right.\\
&&\left. +e^{-\int_0^\infty q_{\phi(x,s)}(\rho_s)ds}e^{\int_0^\infty c(\phi(x,s),\rho_s)ds} \right\},~\forall~x\in S.
\end{eqnarray*}
\item[(b)] The mapping
\begin{eqnarray*}
\rho\in {\cal R}&\rightarrow& W(x,\rho):=\int_0^\infty e^{-\int_0^\tau (q_{\phi(x,s)}(\rho_s)-c(\phi(x,s),\rho_s))ds} \left(\int_S V^\ast(y)\tilde{q}(dy|\phi(x,\tau),\rho_\tau)\right)d\tau\\
&&  +e^{-\int_0^\infty q_{\phi(x,s)}(\rho_s)ds}e^{\int_0^\infty c(\phi(x,s),\rho_s)ds}
\end{eqnarray*}
is lower semicontinuous for each $x\in S.$
\end{itemize}
\end{lemma}
\par\noindent\textit{Proof.}
One can legitimately consider the following DTMDP (discrete-time Markov decision process): according to Lemma 2.29 of \cite{Costa:2013}, all the involved mappings are measurable.
\begin{itemize}
\item The state space is $\textbf{X}:=((0,\infty)\times S)\bigcup\{(\infty,x_\infty)\}$. Whenever the topology is concerned, $(\infty,x_\infty)$ is regarded as an isolated point in $\textbf{X}.$
\item The action space is $\textbf{A}:={\cal R}$.
\item The transition kernel $p$ on ${\cal B}(\textbf{X})$ from $\textbf{X}\times \textbf{A}$, c.f. (\ref{GGZyExponential21}), is given for each $\rho\in \textbf{A}$ by
    \begin{eqnarray*}
    p(\Gamma_1\times\Gamma_2|(\theta,x),\rho)&:=&\int_{\Gamma_2} e^{-\int_0^t q_{\phi(x,s)}(\rho_s)ds}\tilde{q}(\Gamma_1|\phi(x,t),\rho_t)dt,\nonumber\\
    &&~\forall~\Gamma_1\in {\cal B}(S),~\Gamma_2\in {\cal B}((0,\infty)),~x\in S,~\theta\in (0,\infty),\nonumber\\
     p(\{(\infty,x_\infty)\}|(\theta,x),\rho)&:=&e^{-\int_0^\infty q_{\phi(x,s)}(\rho_s)ds},~\forall~x\in S,~\theta\in (0,\infty);\nonumber\\
     p(\{(\infty,x_\infty)\}|(\infty,x_\infty),\rho)&:=&1.\nonumber\\
\end{eqnarray*}
\item The cost function $l$ is a $[0,\infty]$-valued measurable function on $\textbf{X}\times\textbf{A}\times \textbf{X}$ given by
\begin{eqnarray*}\label{ZyExponential55}
l((\theta,x),\rho,(\tau,y)):=\int_0^\infty I\{s<\tau\} c(\phi(x,s),\rho_s)ds,~\forall~((\theta,x),\rho,(\tau,y))\in \textbf{X}\times\textbf{A}\times \textbf{X}.
\end{eqnarray*}
\end{itemize}
The relevant facts and statements for the DTMDP are included in the Appendix.

One can show that under Conditions \ref{GZyExponentialCondition02} and \ref{GGZyExponentialConditionExtra}, for each $(\theta,x)\in \textbf{X}$, $a\in \textbf{A}\rightarrow \int_{\textbf{X}}f(z)p(dz|(\theta,x),a)$ is continuous for each bounded measurable function $f$ on $\textbf{X}$; for each $(\theta,x)\in \textbf{X}$ and $(\tau,y)\in\textbf{X}$, $a\in\textbf{A}\rightarrow l((\theta,x),\rho,(\tau,y))$ is lower semicontinuous, and $\textbf{A}$ is a compact Borel space. Hence, Condition \ref{GGZyExponentialDTMDPCon2} for the DTMDP model $\{\textbf{X},\textbf{A},p,l\}$ is satisfied.

The controlled process in the above DTMDP model  $\{\textbf{X},\textbf{A},p,l\}$  is denoted by $\{Y_n,n=0,1,\dots\}$, where $Y_n=(\Theta_n,X_n)$, and the controlling process is denoted by $\{A_n,n=0,1,\dots\}.$ For $n\ge 1,$ $\Theta_n$ and $X_n$ correspond to the $n$th sojourn time and the post-jump state in the PDMDP, $\Theta_0$ is fictitious, and $X_0$ is the initial state in the PDMDP.
Let $\Sigma$ be the class of all strategies for the DTMDP model $\{\textbf{X},\textbf{A},p,l\}$, and $\Sigma_{DM}^0$ be the class of deterministic Markov strategies in the form $\sigma=(\varphi_n)$ where $\varphi_0((\theta,x))$ does not depend on $\theta\in (0,\infty)$ for each $x\in S.$ We preserve the term of policy for the PDMDP and the term of strategy for the DTMDP.

According to Proposition \ref{GGZyExponentialProposition02}, the function
\begin{eqnarray*}(\theta,x)\in\textbf{X}\rightarrow \textbf{V}^\ast((\theta,x)):=\inf_{\sigma\in \Sigma}\textbf{E}_{(\theta,x)}^{\sigma}\left[e^{\sum_{n=0}^\infty l(Y_n,A_n,Y_{n+1})}\right]
\end{eqnarray*}
is the minimal $[1,\infty]$-valued measurable solution to the optimality equation
\begin{eqnarray*}
\textbf{V}^\ast((\theta,x))&=& \inf_{\rho\in {\cal R}}\left\{\int_0^\infty e^{-\int_0^\tau (q_{\phi(x,s)}(\rho_s)-c(\phi(x,s),\rho_s))ds} \left(\int_S \textbf{V}^\ast((\tau,y))\tilde{q}(dy|\phi(x,\tau),\rho_\tau)\right)d\tau\right.\\
&&\left. +e^{-\int_0^\infty q_{\phi(x,s)}(\rho_s)ds}e^{\int_0^\infty c(\phi(x,s),\rho_s)ds} \right\}
\end{eqnarray*}
for each $x\in S$ and $\theta\in (0,\infty);$ this is just (\ref{ZyExponential02}). Furthermore, by Proposition \ref{GGZyExponentialProposition02}, there exists a deterministic stationary strategy $\sigma^\ast$ for the DTMDP such that $\sigma^\ast((\theta,x))$ attains the above infimum for each $x\in S$ and $\theta\in (0,\infty),$ and any such strategy  $\sigma^\ast$ verifies
\begin{eqnarray*}
\textbf{E}_{(\theta,x)}^{\sigma^\ast}\left[e^{\sum_{n=0}^\infty l(Y_n,A_n,Y_{n+1})}\right]= \inf_{\sigma\in \Sigma}\textbf{E}_{(\theta,x)}^{\sigma}\left[e^{\sum_{n=0}^\infty l(Y_n,A_n,Y_{n+1})}\right],~\forall~(\theta,x)\in\textbf{X}.
\end{eqnarray*}
Let $\hat{\theta}\in~(0,\infty)$ be arbitrarily fixed. The function $\textbf{V}^\ast((\theta,x))$ being measurable in $(\theta,x)\in\textbf{X}$, it follows that $x\in S\rightarrow \textbf{V}^\ast((\hat{\theta},x))$ is measurable.  The strategy $\sigma^\ast$ and the constant $\hat{\theta}$ induce a deterministic Markov strategy $\sigma^{\ast\ast}=(\varphi_n)\in \Sigma^0_{DM}$, where $\varphi_0((\theta,x))=:\sigma^\ast((\hat{\theta},x))$ for each $\theta\in(0,\infty),~x\in S$, and $\varphi_n((\theta,x)):=\sigma((\theta,x))$ for each $n\ge 1$, $\theta\in (0,\infty),~x\in S.$ (The control on the isolated point $(0,x_\infty)$ is irrelevant and we do not specify the definition of the strategy on that point.) This strategy can be identified with a policy $\pi^\ast$ in the PDMDP, c.f. (\ref{GGZZPDMDPStarNumber}). On the other hand, each policy $\pi=(\pi_n)$ can be identified with a deterministic strategy in this DTMDP. Thus,
\begin{eqnarray*}
V^\ast(x)\ge\textbf{V}^\ast((\hat{\theta},x))=\textbf{E}_{(\hat{\theta},x)}^{\sigma^\ast}\left[e^{\sum_{n=0}^\infty l(Y_n,A_n,Y_{n+1})}\right]= \textbf{E}_{(\hat{\theta},x)}^{\sigma^{\ast\ast}}\left[e^{\sum_{n=0}^\infty l(Y_n,A_n,Y_{n+1})}\right]=V(x,\pi^\ast)\ge V^\ast(x)
\end{eqnarray*}
for each $x\in S.$
Consequently, the policy $\pi^\ast$ is optimal, $V^\ast(x)=\textbf{V}^\ast((\theta,x))$ for each $x\in  S$ and $\theta\in(0,\infty);$ recall that $\hat{\theta}$ was arbitrarily fixed. The statement of this lemma now follows. $\hfill\Box$
\bigskip

The policy $\pi^\ast$ in the proof of the previous lemma is actually optimal for problem (\ref{GZPDMDP17}). However, it is not necessarily a deterministic nor stationary policy. Also the reduction of the risk-sensitive PDMDP problem (\ref{GZPDMDP17}) to a risk-sensitive problem for the DTMDP model $\{\textbf{X},\textbf{A},p,l\}$ as seen in the proof of the above theorem will be used without special reference in what follows.

\begin{lemma}\label{GGZyLemmaNew2}
Suppose Conditions \ref{GZyExponentialCondition02}, \ref{GGZyExponentialConditionExtra} and \ref{GGZyExponentialFinitenessCon} are satisfied. For each $x\in S$ and $\rho\in {\cal R}$,
\begin{eqnarray*}
t\in[0,\infty)&\rightarrow& \int_0^t e^{-\int_0^\tau (q_{\phi(x,s)}(\rho_s)-c(\phi(x,s),\rho_s))ds}\int_{S}V^\ast(y)\tilde{q}(dy|\phi(x,\tau),\rho_\tau)d\tau\\
&& +e^{-\int_0^t (q_{\phi(x,s)}(\rho_s)-c(\phi(x,s),\rho_s))ds}V^\ast(\phi(x,t))
\end{eqnarray*}
is monotone nondecreasing in $t\in[0,\infty)$.
\end{lemma}
\par\noindent\textit{Proof.}
 Let $0\le t_1<t_2<\infty$ be arbitrarily fixed.  We need show
\begin{eqnarray}\label{GGZyExponential012}
&&\int_0^{t_2} e^{-\int_0^\tau (q_{\phi(x,s)}(\rho_s)-c(\phi(x,s),\rho_s))ds}\int_{S}V^\ast(y)\tilde{q}(dy|\phi(x,\tau),\rho_\tau)d\tau\nonumber\\
&& +e^{-\int_0^{t_2} (q_{\phi(x,s)}(\rho_s)-c(\phi(x,s),\rho_s))ds}V^\ast(\phi(x,t_2))\nonumber\\
&\ge& \int_0^{t_1} e^{-\int_0^\tau (q_{\phi(x,s)}(\rho_s)-c(\phi(x,s),\rho_s))ds}\int_{S}V^\ast(y)\tilde{q}(dy|\phi(x,\tau),\rho_\tau)d\tau\nonumber\\
&& +e^{-\int_0^{t_1} (q_{\phi(x,s)}(\rho_s)-c(\phi(x,s),\rho_s))ds}V^\ast(\phi(x,t_1)).
\end{eqnarray}
It is without loss of generality to assume
\begin{eqnarray*}
\int_0^{t_2} e^{-\int_0^\tau (q_{\phi(x,s)}(\rho_s)-c(\phi(x,s),\rho_s))ds}\int_{S}V^\ast(y)\tilde{q}(dy|\phi(x,\tau),\rho_\tau)d\tau<\infty.
\end{eqnarray*}
Then all the four terms in (\ref{GGZyExponential012}) are nonnegative and finite, and (\ref{GGZyExponential012}) is equivalent to
\begin{eqnarray}\label{GGZyExponential019}
&&\int_0^{t_2} e^{-\int_0^\tau (q_{\phi(x,s)}(\rho_s)-c(\phi(x,s),\rho_s))ds}\int_{S}V^\ast(y)\tilde{q}(dy|\phi(x,\tau),\rho_\tau)d\tau\nonumber\\
&& +e^{-\int_0^{t_2} (q_{\phi(x,s)}(\rho_s)-c(\phi(x,s),\rho_s))ds}V^\ast(\phi(x,t_2))\nonumber\\
&&-\int_0^{t_1} e^{-\int_0^\tau (q_{\phi(x,s)}(\rho_s)-c(\phi(x,s),\rho_s))ds}\int_{S}V^\ast(y)\tilde{q}(dy|\phi(x,\tau),\rho_\tau)d\tau \nonumber\\
&&-e^{-\int_0^{t_1} (q_{\phi(x,s)}(\rho_s)-c(\phi(x,s),\rho_s))ds}V^\ast(\phi(x,t_1))\nonumber\\
&=&\int_{t_1}^{t_2}e^{-\int_0^\tau (q_{\phi(x,s)}(\rho_s)-c(\phi(x,s),\rho_s))ds}\int_{S}V^\ast(y)\tilde{q}(dy|\phi(x,\tau),\rho_\tau)d\tau\nonumber\\
&&+e^{-\int_0^{t_1} (q_{\phi(x,s)}(\rho_s)-c(\phi(x,s),\rho_s))ds}\left(e^{-\int_{t_1}^{t_2} (q_{\phi(x,s)}(\rho_s)-c(\phi(x,s),\rho_s))ds}V^\ast(\phi(x,t_2))-V^\ast(\phi(x,t_1))\right)\nonumber\\
&=&\left\{\int_0^{t_2-t_1} e^{-\int_0^\tau(q_{\phi(x,s+t_1)}(\rho_{s+t_1})-c(\phi(x,s+t_1),\rho_{s+t_1}))ds}\int_S V^\ast(y)\tilde{q}(dy|\phi(x,t_1+\tau),\rho_{t_1+\tau})d\tau\right.\nonumber\\
&&\left.+e^{-\int_{t_1}^{t_2} (q_{\phi(x,s)}(\rho_s)-c(\phi(x,s),\rho_s))ds}V^\ast(\phi(x,t_2))-V^\ast(\phi(x,t_1))\right\}e^{-\int_0^{t_1} (q_{\phi(x,s)}(\rho_s)-c(\phi(x,s),\rho_s))ds}\nonumber\\
&\ge& 0,
\end{eqnarray}
which is verified as follows.
Let $\delta>0$ be arbitrarily fixed. By Lemma \ref{GGZyPDMDPLem01}, there exists some $\hat{\nu}\in{\cal R}$ such that
\begin{eqnarray*}
V^\ast(\phi(x,t_2))+\delta&\ge& \int_0^\infty \int_S V^\ast(y)\tilde{q}(dy|\phi(x,t_2+\tau),\hat{\nu}_\tau)e^{-\int_0^\tau (q_{\phi(x,t_2+s)}(\hat{\nu}_s)-c(\phi(x,t_2+s),\hat{\nu}_s))ds}d\tau\nonumber\\
&&+e^{-\int_0^\infty q_{\phi(x,t_2+s)}(\hat{\nu}_s)ds}e^{\int_0^\infty c(\phi(x,t_2+s),\hat{\nu}_s)ds}.
\end{eqnarray*}
(Recall $\phi(x,t_2+t)=\phi(\phi(x,t_2),t)$ for each $t \ge 0.$) Consider $\tilde{\nu}\in {\cal R}$ defined by
\begin{eqnarray*}
\tilde{\nu}_s=\left\{\begin{array}{ll}
\rho_{t_1+s}, & \mbox{ if }s\le t_2-t_1;  \\
\hat{\nu}_{s-(t_2-t_1)} & \mbox{ if } s>t_2-t_1. \end{array}\right.
\end{eqnarray*}
Then routine calculations lead to
\begin{eqnarray*}
&&V^\ast(\phi(x,t_1))\\
&\le& \int_0^{t_2-t_1} e^{-\int_0^\tau (q_{\phi(x,t_1+s})(\tilde{\nu}_s)-c(\phi(x,t_1+s),\tilde{\nu}_s))ds} \left(\int_S V^\ast(y)\tilde{q}(dy|\phi(x,t_1+\tau),\tilde{\nu}_\tau)\right)d\tau\\
&&+\int_{t_2-t_1}^\infty e^{-\int_0^\tau (q_{\phi(x,t_1+s})(\tilde{\nu}_s)-c(\phi(x,t_1+s),\tilde{\nu}_s))ds} \left(\int_S V^\ast(y)\tilde{q}(dy|\phi(x,t_1+\tau),\tilde{\nu}_\tau)\right)d\tau\\
&& +e^{-\int_0^{t_2-t_1} (q_{\phi(x,t_1+s)}(\tilde{\nu}_s)-c(\phi(x,t_1+s),\tilde{\nu}_s))ds} e^{-\int_{t_2-t_1}^\infty q_{\phi(x,t_1+s)}(\tilde{\nu}_s)ds}e^{\int_{t_2-t_1}^\infty c(\phi(x,t_1+s),\tilde{\nu}_s)ds}\\
&=&\int_0^{t_2-t_1}e^{-\int_0^\tau (q_{\phi(x,t_1+s)}(\rho_{s+t_1})-c(\phi(x,t_1+s),\rho_{s+t_1}))ds}\int_S V^\ast(y)\tilde{q}(dy|\phi(x,t_1+\tau),\rho_{t_1+\tau})d\tau\\
&&+e^{-\int_0^{t_2-t_1}    (q_{\phi(x,t_1+s)}(\rho_{s+t_1})-c(\phi(x,t_1+s),\rho_{s+t_1}))ds}\\
  &&\times \left\{\int_0^\infty e^{-\int_0^\tau(q_{\phi(x,t_2+s)}(\hat{\nu}_s)-c(\phi(x,t_2+s),\hat{\nu}_s))ds}\int_S V^\ast(y)\tilde{q}(dy|\phi(x,t_2+\tau),\hat{\nu}_\tau)d\tau\right.\\
  &&\left.
  +e^{-\int_0^\infty q_{\phi(x,t_2+s)}(\hat{\nu}_s)ds}e^{\int_0^\infty c(\phi(x,t_2+s),\hat{\nu}_s)ds}\right\}\\
&\le&\int_0^{t_2-t_1}e^{-\int_0^\tau (q_{\phi(x,t_1+s)}(\rho_{s+t_1})-c(\phi(x,t_1+s),\rho_{s+t_1}))ds}\int_S V^\ast(y)\tilde{q}(dy|\phi(x,t_1+\tau),\rho_{t_1+\tau})d\tau\\
&&+e^{-\int_0^{t_2-t_1}    (q_{\phi(x,t_1+s)}(\rho_{s+t_1})-c(\phi(x,t_1+s),\rho_{s+t_1}))ds} (V^\ast(\phi(x,t_2))+\delta).
\end{eqnarray*}
Since $\delta>0$ was arbitrarily fixed, now it follows that the term in the parenthesis in (\ref{GGZyExponential019}) is nonnegative, and thus inequality (\ref{GGZyExponential019}) is verified. $\hfill\Box$
\bigskip

\begin{lemma}\label{GGZyExponentialLemma05}
Suppose Conditions \ref{GZyExponentialCondition02}, \ref{GGZyExponentialConditionExtra} and \ref{GGZyExponentialFinitenessCon} are satisfied. For each $x\in S$, there is some $\rho^\ast\in {\cal R}$ such that
\begin{eqnarray}\label{GGZyExponential02}
V^\ast(x)&=&\inf_{\rho\in{\cal  R}}\left\{\int_0^t e^{-\int_0^s(q_{\phi(x,v)}(\rho_v)-c(\phi(x,v),\rho_v))dv}\int_S V^\ast(y)\tilde{q}(dy|\phi(x,s),\rho_s)ds\right.\nonumber\\
&&\left.+e^{-\int_0^t (q_{\phi(x,s)}(\rho_s)-c(\phi(x,s),\rho_s))ds}V^\ast(\phi(x,t))\right\}\nonumber\\
&=&\int_0^t e^{-\int_0^s(q_{\phi(x,v)}(\rho^\ast_v)-c(\phi(x,v),\rho^\ast_v))dv}\int_S V^\ast(y)\tilde{q}(dy|\phi(x,s),\rho^\ast_s)ds\nonumber\\
&&+e^{-\int_0^t (q_{\phi(x,s)}(\rho^\ast_s)-c(\phi(x,s),\rho^\ast_s))ds}V^\ast(\phi(x,t)),~\forall~t\ge 0.
\end{eqnarray}
\end{lemma}
\par\noindent\textit{Proof.} Let $x\in S$ be fixed, and let $\rho^\ast\in {\cal R}$ be such that $V^\ast(x)=W(x,\rho^\ast)$, see Lemma \ref{GGZyPDMDPLem01}. Suppose $t\in[0,\infty)$ is arbitrarily fixed. Consider $\tilde{\rho}\in {\cal R}$ defined by
$
\tilde{\rho}_s=\rho^\ast_{t+s}$ for each $s>0$.
Then
\begin{eqnarray*}
V^\ast(x)&=& \int_0^t e^{-\int_0^s(q_{\phi(x,v)}(\rho^\ast_v)-c(\phi(x,v),\rho^\ast_v))dv}\int_S V^\ast(y)\tilde{q}(dy|\phi(x,s),\rho^\ast_s)ds+e^{-\int_0^t (q_{\phi(x,s)}(\rho^\ast_s)-c(\phi(x,s),\rho^\ast_s))ds}\\
&&\times\left\{\int_0^\infty e^{-\int_0^\tau (q_{\phi(x,t+s)}(\tilde{\rho}_s)-c(\phi(x,s+t),\tilde{\rho}_s))ds}\int_S V^\ast(y)\tilde{q}(dy|\phi(x,\tau+t),\tilde{\rho}_\tau)d\tau\right.\\
&&\left.+e^{-\int_0^\infty q_{\phi(x,t+s)}(\tilde{\rho}_s)ds}e^{-\int_0^\infty c(\phi(x,t+s),\tilde{\rho}_s)ds}\right\}\\
&\ge&\int_0^t e^{-\int_0^s(q_{\phi(x,v)}(\rho^\ast_v)-c(\phi(x,v),\rho^\ast_v))dv}\int_S V^\ast(y)\tilde{q}(dy|\phi(x,s),\rho^\ast_s)ds\\
&&+e^{-\int_0^t (q_{\phi(x,s)}(\rho^\ast_s)-c(\phi(x,s),\rho^\ast_s))ds}V^\ast(\phi(x,t));
\end{eqnarray*}
recall (\ref{GGZyExponential01}).
On the other hand, by Lemma \ref{GGZyLemmaNew2},
\begin{eqnarray*}
V^\ast(x)&\le&\inf_{\rho\in{\cal  R}}\left\{\int_0^t e^{-\int_0^s(q_{\phi(x,v)}(\rho_v)-c(\phi(x,v),\rho_v))dv}\int_S V^\ast(y)\tilde{q}(dy|\phi(x,s),\rho_s)ds\right.\\
&&\left.+e^{-\int_0^t (q_{\phi(x,s)}(\rho_s)-c(\phi(x,s),\rho_s))ds}V^\ast(\phi(x,t))\right\}.
\end{eqnarray*}
The statement of this lemma is thus proved. $\hfill\Box$
\bigskip

 \begin{lemma}\label{GGZyExponentialLemma06}
 Suppose Conditions \ref{GZyExponentialCondition02}, \ref{GGZyExponentialConditionExtra} and \ref{GGZyExponentialFinitenessCon} are satisfied. Then for each $x\in S,$ $t\in[0,\infty)\rightarrow V^\ast(\phi(x,t))$ is absolutely continuous.
 \end{lemma}
 \par\noindent\textit{Proof.} This immediately follows from Lemma \ref{GGZyExponentialLemma05}. $\hfill\Box$
 \bigskip

\par\noindent\textit{Proof of Theorem \ref{GGZyExponentialTheorem}.}
(a) Under Conditions \ref{GZyExponentialCondition02}, \ref{GGZyExponentialConditionExtra} and \ref{GGZyExponentialFinitenessCon}, by Lemma \ref{GGZyExponentialLemma06}, for each $x\in S,$ let $t\in[0,\infty)\rightarrow U^\ast(x,t)$ be an integrable real-valued function such that $U^\ast(x,t)$ coincides with the derivative of $t\in[0,\infty)\rightarrow V(\phi(x,t))$ almost everywhere. Let $x\in S$ and $t\in[0,\infty)$ be fixed, and let $\rho^\ast\in {\cal R}$ be from Lemma \ref{GGZyExponentialLemma05}.

By Lemmas \ref{GGZyExponentialLemma05} and \ref{GGZyExponentialLemma06},
\begin{eqnarray*}
\int_0^\tau e^{-\int_0^s(q_{\phi(x,v)}(\rho^\ast_v)-c(\phi(x,v),\rho^\ast_v))dv}\int_S V^\ast(y)\tilde{q}(dy|\phi(x,s),\rho^\ast_s)ds
\end{eqnarray*}
and
\begin{eqnarray*}
e^{-\int_0^\tau (q_{\phi(x,s)}(\rho^\ast_s)-c(\phi(x,s),\rho^\ast_s))ds}V^\ast(\phi(x,\tau))
\end{eqnarray*}
are absolutely continuous in $\tau$ and are finite for each $\tau\in[0,\infty)$.
Since $\phi(x,0)=x$, see (\ref{GGZyExponential01}),
\begin{eqnarray*}
&&e^{-\int_0^t (q_{\phi(x,s)}(\rho^\ast_s)-c(\phi(x,s),\rho^\ast_s))ds}V^\ast(\phi(x,t))-V^\ast(x)\nonumber\\
&=&\int_0^t e^{-\int_0^\tau (q_{\phi(x,s)}(\rho^\ast_s)-c(\phi(x,s),\rho^\ast_s) )ds}\left\{U^\ast(x,\tau)- (q_{\phi(x,\tau)}(\rho^\ast_\tau)-c(\phi(x,\tau),\rho^\ast_\tau) )V^\ast(\phi(x,\tau))     \right\}d\tau.
\end{eqnarray*}
Now by Lemma \ref{GGZyExponentialLemma05},
\begin{eqnarray}\label{GGZyExponential08}
0
&=&\int_0^t e^{-\int_0^s(q_{\phi(x,v)}(\rho^\ast_v)-c(\phi(x,v),\rho^\ast_v))dv}\int_S V^\ast(y)\tilde{q}(dy|\phi(x,s),\rho^\ast_s)ds\nonumber\\
&&+e^{-\int_0^t (q_{\phi(x,s)}(\rho^\ast_s)-c(\phi(x,s),\rho^\ast_s))ds}V^\ast(\phi(x,t))-V^\ast(x)\nonumber\\
&=&\int_0^t e^{-\int_0^\tau(q_{\phi(x,v)}(\rho^\ast_v)-c(\phi(x,v),\rho^\ast_v))dv}\left\{\int_S V^\ast(y)\tilde{q}(dy|\phi(x,\tau),\rho^\ast_\tau)+U^\ast(x,\tau)\right.\nonumber\\
&&\left.- (q_{\phi(x,\tau)}(\rho^\ast_\tau)-c(\phi(x,\tau),\rho^\ast_\tau) )V^\ast(\phi(x,\tau)) \right\}d\tau \nonumber\\
&\ge& \int_0^t e^{-\int_0^\tau(q_{\phi(x,v)}(\rho^\ast_v)-c(\phi(x,v),\rho^\ast_v))dv}\left\{U^\ast(x,\tau)\right. \nonumber\\
&&+\left.\inf_{a\in A}\left\{ \int_S V^\ast(y)\tilde{q}(dy|\phi(x,\tau),a)- (q_{\phi(x,\tau)}(a)-c(\phi(x,\tau),a) )V^\ast(\phi(x,\tau))\right\} \right\}d\tau\nonumber\\
&=& \int_0^t e^{-\int_0^\tau(q_{\phi(x,v)}(\rho^\ast_v)-c(\phi(x,v),\rho^\ast_v))dv}\left\{U^\ast(x,\tau)+\int_S V^\ast(y)\tilde{q}(dy|\phi(x,\tau),f(\phi(x,\tau)))\right.\nonumber\\
&&\left. - (q_{\phi(x,\tau)}(f(\phi(x,\tau)))-c(\phi(x,\tau),f(\phi(x,\tau))) )V^\ast(\phi(x,\tau)) \right\}d\tau,
\end{eqnarray}
where $f$ is a measurable mapping from $S$ to $A$ such that
\begin{eqnarray*}
&&\inf_{a\in A}\left\{ \int_S V^\ast(y)\tilde{q}(dy|x,a)- (q_{x}(a)-c(x,a) )V^\ast(x)\right\}\\
&=& \int_S V^\ast(y)\tilde{q}(dy|x,f(x))- (q_{x}(\varphi(x))-c(x,f(x)) )V^\ast(x)
\end{eqnarray*}
for each $x\in S$; the existence of such a mapping is according to a well known measurable selection theorem, c.f. Proposition D.5 of \cite{Hernandez-Lerma:1996}.

Note that $e^{-\int_0^\tau(q_{\phi(x,v)}(\rho_v)-c(\phi(x,v),\rho_v))dv}$ is bounded and separated from zero in $\tau\in[0,t]$ for each $\rho\in {\cal R};$ recall Condition \ref{GGZyExponentialConditionExtra}. So
\begin{eqnarray*}\int_0^t e^{-\int_0^\tau(q_{\phi(x,v)}(\rho^\ast_v)-c(\phi(x,v),\rho^\ast_v))dv}\left\{U^\ast(x,\tau) - (q_{\phi(x,\tau)}(f(\phi(x,\tau)))-c(\phi(x,\tau),f(\phi(x,\tau))) )V^\ast(\phi(x,\tau)) \right\}d\tau
\end{eqnarray*}
is finite.
If \begin{eqnarray*}
\int_0^t  \int_S V^\ast(y)\tilde{q}(dy|\phi(x,\tau),f(\phi(x,\tau)))d\tau=\infty,
\end{eqnarray*}
then
\begin{eqnarray*}&&\int_0^t e^{-\int_0^\tau(q_{\phi(x,v)}(\rho^\ast_v)-c(\phi(x,v),\rho^\ast_v))dv}\left\{U^\ast(x,\tau)+\int_S V^\ast(y)\tilde{q}(dy|\phi(x,\tau),f(\phi(x,\tau)))\right.\\
&&\left. - (q_{\phi(x,\tau)}(f(\phi(x,\tau)))-c(\phi(x,\tau),f(\phi(x,\tau))) )V^\ast(\phi(x,\tau)) \right\}d\tau=\infty,
\end{eqnarray*}
which is against (\ref{GGZyExponential08}). Therefore,
\begin{eqnarray*}
\int_0^t  \int_S V^\ast(y)\tilde{q}(dy|\phi(x,\tau),f(\phi(x,\tau)))d\tau<\infty.
\end{eqnarray*}
Then
\begin{eqnarray*}
&&\int_0^v e^{-\int_0^\tau (q_{\phi(x,s)}(f(\phi(x,s)))-c(\phi(x,s),f(\phi(x,s))))ds}\int_{S}V^\ast(y)\tilde{q}(dy|\phi(x,\tau), f(\phi(x,\tau)))d\tau\\
&& +e^{-\int_0^v (q_{\phi(x,s)}(f(\phi(x,s)))-c(\phi(x,s),f(\phi(x,s))))ds}V^\ast(\phi(x,v))
\end{eqnarray*}
is absolutely continuous on $[0,t].$ After legitimately differentiating the above expression with respect to $v$, and applying Lemma \ref{GGZyLemmaNew2}, we see
\begin{eqnarray*}
&&U^\ast(x,v)+\int_S V^\ast(y)\tilde{q}(dy|\phi(x,v),f(\phi(x,v)))\nonumber\\
&&- (q_{\phi(x,v)}(f(\phi(x,v)))-c(\phi(x,v),f(\phi(x,v))) )V^\ast(\phi(x,v))\ge 0
\end{eqnarray*}
for almost all $v\in[0,t].$ This and (\ref{GGZyExponential08}) imply
\begin{eqnarray*}
U^\ast(x,\tau)+\inf_{a\in A}\left\{ \int_S V^\ast(y)\tilde{q}(dy|\phi(x,\tau),a)- (q_{\phi(x,\tau)}(a)-c(\phi(x,\tau),a) )V^\ast(\phi(x,\tau))\right\}=0
\end{eqnarray*}
almost everywhere in $\tau\in[0,t].$ Remember, $t\in[0,\infty)$ was arbitrarily fixed. The first part of (a) is thus verified, and we postpone the justification of the second part of (a) after the proof of part (b).

(b) We use the same notation as in the above.
Note that \begin{eqnarray}\label{GGZyExponential09}
\liminf_{t\rightarrow \infty}\left\{e^{-\int_0^t (q_{\phi(x,s)}(f(\phi(x,s)))- c(\phi(x,s),f(\phi(x,s))))ds}\right\}\ge e^{-\int_0^\infty q_{\phi(x,s)}(f(\phi(x,s)))ds} e^{\int_0^\infty c(\phi(x,s),f(\phi(x,s))))ds}.
\end{eqnarray}
Indeed, if either $\int_0^\infty q_{\phi(x,s)}(f(\phi(x,s)))ds$ or $\int_0^\infty c(\phi(x,s),f(\phi(x,s))))ds$ is finite, then in the above inequality, the equality takes place; and if both  $\int_0^\infty q_{\phi(x,s)}(f(\phi(x,s)))ds$ and $\int_0^\infty c(\phi(x,s),f(\phi(x,s))))ds$ are infinite, then the right hand side of the inequality is zero according to (\ref{GGZyExponential56}).

In the proof of part (a), it was observed that
\begin{eqnarray*}
\int_0^t e^{-\int_0^s(q_{\phi(x,v)}(f(\phi(x,v)))-c(\phi(x,v),f(\phi(x,v))))dv}\int_S V^\ast(y)\tilde{q}(dy|\phi(x,s),f(\phi(x,s)))ds
\end{eqnarray*}
and
\begin{eqnarray*}
e^{-\int_0^t (q_{\phi(x,s)}(f(\phi(x,s)))-c(\phi(x,s),f(\phi(x,s))))ds}V^\ast(\phi(x,t))
\end{eqnarray*}
are absolutely continuous in $t$ and are thus finite for each $t\in[0,\infty)$. As in the proof of part (a), similar calculations to those in (\ref{GGZyExponential08}) imply that for each $t\in[0,\infty),$
\begin{eqnarray*}
&&\int_0^t e^{-\int_0^s(q_{\phi(x,v)}(f(\phi(x,v)))-c(\phi(x,v),f(\phi(x,v))))dv}\int_S V^\ast(y)\tilde{q}(dy|\phi(x,s),f(\phi(x,s)))ds\nonumber\\
&&+e^{-\int_0^t (q_{\phi(x,s)}(f(\phi(x,s)))-c(\phi(x,s),f(\phi(x,s))))ds}V^\ast(\phi(x,t))-V^\ast(x)\nonumber\\
&=& \int_0^t e^{-\int_0^\tau(q_{\phi(x,v)}(f(\phi(x,v)))-c(\phi(x,v),f(\phi(x,v))))dv}\left\{U^\ast(x,\tau)+\int_S V^\ast(y)\tilde{q}(dy|\phi(x,\tau),f(\phi(x,\tau)))\right.\nonumber\\
&&\left. - (q_{\phi(x,\tau)}(f(\phi(x,\tau)))-c(\phi(x,\tau),f(\phi(x,\tau)) ))V^\ast(\phi(x,\tau)) \right\}d\tau=0,
\end{eqnarray*}
where the last equality is by what was established in part (a). Therefore, for each $t\in[0,\infty),$
\begin{eqnarray*}
&&V^\ast(x)-\int_0^t e^{-\int_0^s(q_{\phi(x,v)}(f(\phi(x,v)))-c(\phi(x,v),f(\phi(x,v))))dv}\int_S V^\ast(y)\tilde{q}(dy|\phi(x,s),f(\phi(x,s)))ds\\
&=&e^{-\int_0^t (q_{\phi(x,s)}(f(\phi(x,s)))-c(\phi(x,s),f(\phi(x,s))))ds}V^\ast(\phi(x,t))\\
&\ge &e^{-\int_0^t (q_{\phi(x,s)}(f(\phi(x,s)))-c(\phi(x,s),f(\phi(x,s))))ds},
\end{eqnarray*}
where the inequality holds because $V^\ast(x)\ge 1$ for each $x\in S.$ Taking $\liminf_{t\rightarrow \infty}$ on the both sides of the previous equality yields:
\begin{eqnarray*}
&&V^\ast(x)-\int_0^\infty e^{-\int_0^s(q_{\phi(x,v)}(f(\phi(x,v)))-c(\phi(x,v),f(\phi(x,v))))dv}\int_S V^\ast(y)\tilde{q}(dy|\phi(x,s),f(\phi(x,s)))ds\\
&\ge &e^{-\int_0^\infty q_{\phi(x,s)}(f(\phi(x,s)))ds} e^{\int_0^\infty c(\phi(x,s),f(\phi(x,s))))ds}
\end{eqnarray*}
with the inequality following from (\ref{GGZyExponential09}). Hence
\begin{eqnarray*}
V^\ast(x)&\ge&\int_0^\infty e^{-\int_0^s(q_{\phi(x,v)}(f(\phi(x,v)))-c(\phi(x,v),f(\phi(x,v))))dv}\int_S V^\ast(y)\tilde{q}(dy|\phi(x,s),f(\phi(x,s)))ds\\
&&+e^{-\int_0^\infty q_{\phi(x,s)}(f(\phi(x,s)))ds} e^{\int_0^\infty c(\phi(x,s),f(\phi(x,s))))ds}=W(x,\tilde{f}^x)\ge V^\ast(x).
\end{eqnarray*}
Here it is clear that $s\in[0,\infty)\rightarrow f(\phi(x,s))$ can be identified as an element of ${\cal R}$, denoted as $\tilde{f}^x$. In fact, $\tilde{f}_s^x=\delta_{\{f(\phi(x,s))\}}$ for each $s\in [0,\infty)$, whereas $x\in S\rightarrow \tilde{f}^x\in {\cal R}$ is measurable. This measurable mapping $x\in S\rightarrow \tilde{f}^x\in {\cal R}$ defines a deterministic stationary optimal strategy for the risk-sensitive DTMDP problem (\ref{ZyExponential02}) by Proposition \ref{GGZyExponentialProposition02}. It is clear that the measurable mapping $x\in S\rightarrow f(x)\in A$ defines an optimal deterministic stationary policy for the PDMDP problem (\ref{GZPDMDP17}).

Finally, we show the remaining part of (a). Let $H^\ast$ be a measurable $[1,\infty)$-valued function on $S$ such that
 \begin{eqnarray*}
&& -(H^\ast(\phi(x,t))-H^\ast(x))\\
&=&\int_0^t \inf_{a\in A}\left\{ \int_S H^\ast(y)\tilde{q}(dy|\phi(x,\tau),a)- (q_{\phi(x,\tau)}(a)-c(\phi(x,\tau),a) )H^\ast(\phi(x,\tau))\right\}d\tau,\\
&&~t\in[0,\infty),x\in S.
 \end{eqnarray*}
There exists a measurable mapping $h$ from $S$ to $A$ such that
 \begin{eqnarray*}
&&\inf_{a\in A}\left\{ \int_S H^\ast(y)\tilde{q}(dy|x,a)- (q_{x}(a)-c(x,a) )H^\ast(x)\right\}\\
&=&\int_S H^\ast(y)\tilde{q}(dy|x,h(x))- (q_{x}(h(x))-c(x,h(x)) )H^\ast(x),~\forall~x\in S;
 \end{eqnarray*}
c.f., Proposition D.5 of \cite{Hernandez-Lerma:1996}. It follows that $\int_0^s\int_S H^\ast(y)\tilde{q}(dy|\phi(x,\tau),h(\phi(x,\tau)))d\tau$ is absolutely continuous in $s\in[0,t]$ for each $t\ge 0.$ As in the proof of part (b),
\begin{eqnarray*}
&&\int_0^t e^{-\int_0^s(q_{\phi(x,v)}(h(\phi(x,v)))-c(\phi(x,v),h(\phi(x,v))))dv}\int_S H^\ast(y)\tilde{q}(dy|\phi(x,s),h(\phi(x,s)))ds\nonumber\\
&&+e^{-\int_0^t (q_{\phi(x,s)}(h(\phi(x,s)))-c(\phi(x,s),h(\phi(x,s))))ds}H^\ast(\phi(x,t))-H^\ast(x) =0,~\forall~t\in[0,\infty),
\end{eqnarray*}
and by passing to the lower limit as $t\rightarrow \infty$,
\begin{eqnarray}\label{GGZYPDMDP12}
H^\ast(x)&\ge& \int_0^\infty e^{-\int_0^s(q_{\phi(x,v)}(h(\phi(x,v)))-c(\phi(x,v),h(\phi(x,v))))dv}\int_S H^\ast(y)\tilde{q}(dy|\phi(x,s),h(\phi(x,s)))ds\nonumber\\
&&+e^{-\int_0^\infty q_{\phi(x,s)}(h(\phi(x,s)))ds} e^{\int_0^\infty c(\phi(x,s),h(\phi(x,s))))ds}\nonumber\\
&\ge&\inf_{\rho\in {\cal R}}\left\{\int_0^\infty e^{-\int_0^\tau (q_{\phi(x,s)}(\rho_s)-c(\phi(x,s),\rho_s))ds} \left(\int_S H^\ast(y)\tilde{q}(dy|\phi(x,\tau),\rho_\tau)\right)d\tau\right.\nonumber\\
&&\left. +e^{-\int_0^\infty q_{\phi(x,s)}(\rho_s)ds}e^{\int_0^\infty c(\phi(x,s),\rho_s)ds} \right\},~\forall~x\in S.
\end{eqnarray}
It remains to refer to Proposition \ref{GGZyExponentialProposition02} for that $H^\ast(x)\ge V^\ast(x)$ for each $x\in S.$
$\hfill\Box$
\bigskip

\par\noindent\textit{Proof of Theorem \ref{GGGZyExponentialTheorem2}}.
Let $V^\ast_0(x):=1$ for each $x\in S.$ For each $n\ge 0,$ one can legitimately define
\begin{eqnarray}\label{GGZYPDMDP15}
V^\ast_{n+1}(x)&=& \inf_{\rho\in {\cal R}}\left\{\int_0^\infty e^{-\int_0^\tau (q_{\phi(x,s)}(\rho_s)-c(\phi(x,s),\rho_s))ds} \left(\int_S V^\ast_n(y)\tilde{q}(dy|\phi(x,\tau),\rho_\tau)\right)d\tau\right.\nonumber\\
&&\left. +e^{-\int_0^\infty q_{\phi(x,s)}(\rho_s)ds}e^{\int_0^\infty c(\phi(x,s),\rho_s)ds} \right\},~\forall~x\in S.
\end{eqnarray}
Recall that the DTMDP model $\{\textbf{X},\textbf{A},p,l\}$ satisfies Condition \ref{GGZyExponentialDTMDPCon2}, as noted in the proof of Lemma \ref{GGZyPDMDPLem01}.
Then by Proposition \ref{GGZyExponentialProposition02},  $\{V_n^\ast\}$ is a monotone nondecreasing sequence of $[1,\infty)$-valued measurable functions on $S$ such that $V^\ast_n(x)\uparrow V^\ast(x)$ as $n\uparrow \infty,$ for each $x\in S.$

Let $n\ge 0$ be fixed.
As in Lemma \ref{GGZyExponentialLemma05},
for each $x\in S$, there is some $\rho^\ast\in {\cal R}$ such that
\begin{eqnarray*}
V^\ast_{n+1}(x)&=&\inf_{\rho\in{\cal  R}}\left\{\int_0^t e^{-\int_0^s(q_{\phi(x,v)}(\rho_v)-c(\phi(x,v),\rho_v))dv}\int_S V^\ast_n(y)\tilde{q}(dy|\phi(x,s),\rho_s)ds\right.\nonumber\\
&&\left.+e^{-\int_0^t (q_{\phi(x,s)}(\rho_s)-c(\phi(x,s),\rho_s))ds}V^\ast_{n+1}(\phi(x,t))\right\}\nonumber\\
&=&\int_0^t e^{-\int_0^s(q_{\phi(x,v)}(\rho^\ast_v)-c(\phi(x,v),\rho^\ast_v))dv}\int_S V_n^\ast(y)\tilde{q}(dy|\phi(x,s),\rho^\ast_s)ds\nonumber\\
&&+e^{-\int_0^t (q_{\phi(x,s)}(\rho^\ast_s)-c(\phi(x,s),\rho^\ast_s))ds}V_{n+1}^\ast(\phi(x,t)),~\forall~t\ge 0.
\end{eqnarray*}
Also the relevant version of Lemma  \ref{GGZyLemmaNew2} holds: for each $x\in S$ and $\rho\in {\cal R}$,
\begin{eqnarray*}
t\in[0,\infty)&\rightarrow& \int_0^t e^{-\int_0^\tau (q_{\phi(x,s)}(\rho_s)-c(\phi(x,s),\rho_s))ds}\int_{S}V^\ast_n(y)\tilde{q}(dy|\phi(x,\tau),\rho_\tau)d\tau\\
&& +e^{-\int_0^t (q_{\phi(x,s)}(\rho_s)-c(\phi(x,s),\rho_s))ds}V_{n+1}^\ast(\phi(x,t))
\end{eqnarray*}
is monotone nondecreasing in $t\in[0,\infty)$.
Clearly, $V^\ast_{n+1}(\phi(x,t))$ is absolutely continuous in $t\in[0,\infty)$ for each $x\in S$.

Corresponding to (\ref{GGZyExponential08}), we now have
\begin{eqnarray*}
0
&=&\int_0^t e^{-\int_0^s(q_{\phi(x,v)}(\rho^\ast_v)-c(\phi(x,v),\rho^\ast_v))dv}\int_S V^\ast_n(y)\tilde{q}(dy|\phi(x,s),\rho^\ast_s)ds\nonumber\\
&&+e^{-\int_0^t (q_{\phi(x,s)}(\rho^\ast_s)-c(\phi(x,s),\rho^\ast_s))ds}V^\ast_{n+1}(\phi(x,t))-V_{n+1}^\ast(x)\nonumber\\
&=&\int_0^t e^{-\int_0^\tau(q_{\phi(x,v)}(\rho^\ast_v)-c(\phi(x,v),\rho^\ast_v))dv}\left\{\int_S V_n^\ast(y)\tilde{q}(dy|\phi(x,\tau),\rho^\ast_\tau)+U_{n+1}^\ast(x,\tau)\right.\nonumber\\
&&\left.- (q_{\phi(x,\tau)}(\rho^\ast_\tau)-c(\phi(x,\tau),\rho^\ast_\tau) )V_{n+1}^\ast(\phi(x,\tau)) \right\}d\tau \nonumber\\
&\ge& \int_0^t e^{-\int_0^\tau(q_{\phi(x,v)}(\rho^\ast_v)-c(\phi(x,v),\rho^\ast_v))dv}\left\{U_{n+1}^\ast(x,\tau)\right. \nonumber\\
&&+\left.\inf_{a\in A}\left\{ \int_S V_n^\ast(y)\tilde{q}(dy|\phi(x,\tau),a)- (q_{\phi(x,\tau)}(a)-c(\phi(x,\tau),a) )V_{n+1}^\ast(\phi(x,\tau))\right\} \right\}d\tau\nonumber\\
&=& \int_0^t e^{-\int_0^\tau(q_{\phi(x,v)}(\rho^\ast_v)-c(\phi(x,v),\rho^\ast_v))dv}\left\{U_{n+1}^\ast(x,\tau)+\int_S V_n^\ast(y)\tilde{q}(dy|\phi(x,\tau),f(\phi(x,\tau)))\right.\nonumber\\
&&\left. - (q_{\phi(x,\tau)}(f(\phi(x,\tau)))-c(\phi(x,\tau),f(\phi(x,\tau)) )V_{n+1}^\ast(\phi(x,\tau)) \right\}d\tau,
\end{eqnarray*}
where $\tau\in[0,t]\rightarrow U^\ast_{n+1}(x,\tau)$ is integrable and coincides with $\frac{\partial V^\ast_{n+1}(\phi(x,t))}{\partial t}$ almost everywhere, and $f$ is some measurable mapping from $S$ to $A$, whose existence is guaranteed by Proposition D.5 of \cite{Hernandez-Lerma:1996}.
Continued from the above relation, the reasoning in the proof of the first assertion in part (a) of Theorem \ref{GGZyExponentialTheorem} can be followed: eventually we see
\begin{eqnarray*}
U_{n+1}^\ast(x,\tau)+\inf_{a\in A}\left\{ \int_S V_n^\ast(y)\tilde{q}(dy|\phi(x,\tau),a)- (q_{\phi(x,\tau)}(a)-c(\phi(x,\tau),a) )V_{n+1}^\ast(\phi(x,\tau))\right\}=0
\end{eqnarray*}
almost everywhere in $\tau\in[0,t],$ i.e., the equation
\begin{eqnarray}\label{GGPDMDPEqn10}
 && -(V(\phi(x,t))-V(x))\nonumber\\
&=&\int_0^t \inf_{a\in A}\left\{ \int_S V^{\ast}_n(y)\tilde{q}(dy|\phi(x,\tau),a)- (q_{\phi(x,\tau)}(a)-c(\phi(x,\tau),a) )V(\phi(x,\tau))\right\}d\tau,\nonumber\\
&&~t\in[0,\infty),x\in S,
 \end{eqnarray}
is satisfied by $V=V^\ast_{n+1}.$

Recall that $V^\ast_{0}=V^{(0)}$. Suppose the recursive definition in (\ref{GZYPDMDPEqn11}) is valid up to step $n$, and $V^\ast_{n}(x)=V^{(n)}(x)$ for each $x\in S.$  Consider an arbitrarily fixed $[1,\infty)$-valued measurable solution $V$ to (\ref{GGPDMDPEqn10}), and let $f^\ast$ be a measurable mapping from $S$ to $A$ such that
\begin{eqnarray*}
&&\inf_{a\in A}\left\{ \int_S V^\ast_n(y)\tilde{q}(dy|x,a)- (q_{x}(a)-c(x,a) )V(x)\right\}\\
&=&\int_S V_n^\ast(y)\tilde{q}(dy|x,f^\ast(x))- (q_{x}(f^\ast(x))-c(x,f^\ast(x)) )V(x),~\forall~x\in S.
\end{eqnarray*}

One can follow the reasoning in the last part of the proof of Theorem \ref{GGZyExponentialTheorem}, and see, c.f. (\ref{GGZYPDMDP12}),
\begin{eqnarray*}
V(x)&\ge& \int_0^\infty e^{-\int_0^s(q_{\phi(x,v)}(f^\ast(\phi(x,v)))-c(\phi(x,v),f^\ast(\phi(x,v))))dv}\int_S V_n^\ast(y)\tilde{q}(dy|\phi(x,s),f^\ast(\phi(x,s)))ds\\
&&+e^{-\int_0^\infty q_{\phi(x,s)}(f^\ast(\phi(x,s)))ds} e^{\int_0^\infty c(\phi(x,s),f^\ast(\phi(x,s))))ds}\\
&\ge&\inf_{\rho\in {\cal R}}\left\{\int_0^\infty e^{-\int_0^\tau (q_{\phi(x,s)}(\rho_s)-c(\phi(x,s),\rho_s))ds} \left(\int_S V^\ast_n(y)\tilde{q}(dy|\phi(x,\tau),\rho_\tau)\right)d\tau\right.\\
&&\left. +e^{-\int_0^\infty q_{\phi(x,s)}(\rho_s)ds}e^{\int_0^\infty c(\phi(x,s),\rho_s)ds} \right\}=V^\ast_{n+1}(x),~\forall~x\in S,
\end{eqnarray*}
where the last equality is by (\ref{GGZYPDMDP15}).
Thus, $V^\ast_{n+1}$ is the minimal $[1,\infty)$-valued measurable solution to (\ref{GGPDMDPEqn10}), and coincides with $V^{(n+1)}$.  Therefore, by induction $V^\ast_{n}=V^{(n)}$ for each $n\ge 0.$ It follows now that $V^{(n)}(x)\uparrow V^\ast(x)$ as $n\uparrow \infty$ for each $x\in S.$ $\hfill\Box$
\bigskip

\section{Conclusion}\label{GGZPDMDPSecConclusion}
In this paper, we considered total undiscounted risk-sensitive PDMDP in Borel state and action spaces with a nonnegative cost rate. The transition and cost rates are assumed to be locally integrable along the drift. Under quite natural conditions, we showed that the value function is a solution to the optimality equation, justified the value iteration algorithm, and showed the existence of deterministic stationary optimal policy. As a corollary, the obtained results were applied to improving significantly known results for finite horizon undiscounted and infinite horizon discounted risk-sensitive CTMDP in the literature.

\appendix
\section{Appendix}
For ease of reference, we present the relevant notations and facts about the risk-sensitive problem for a DTMDP. The proofs of the presented statements can be found in \cite{Jaskiewicz:2008} or \cite{Zhang:2017}. Standard description of a DTMDP can be found in e.g., \cite{Hernandez-Lerma:1996,Piunovskiy:1997}.

Consider a discrete-time Markov decision process with the following primitives:
\begin{itemize}
\item $\textbf{X}$ is a nonempty Borel state space.
\item $\textbf{A}$ is a nonempty Borel action space.
\item $p(dy|x,a)$ is a stochastic kernel on ${\cal B}(\textbf{X})$ given $(x,a)\in \textbf{X}\times\textbf{A}$.
\item $l$ a $[0,\infty]$-valued measurable cost function on $\textbf{X}\times\textbf{A}\times\textbf{X}.$
\end{itemize}
Let $\Sigma$ be the space of strategies, and $\Sigma_{DM}$ be the space of all deterministic strategies for the DTMDP.
Let the controlled and controlling processes be denoted by $\{Y_n, n=0,1,\dots,\infty\}$ and $\{A_n,n=0,1,\dots,\infty\}$, respectively.
The strategic measure of a strategy $\sigma$ given the initial state $x\in \textbf{X}$ is denoted by $\textbf{P}_x^\sigma$. The expectation taken with respect to $\textbf{P}_x^\sigma$ is denoted by $\textbf{E}_x^\sigma.$

Consider the optimal control problem
\begin{eqnarray}\label{GGZyExponentialProblem2}
\mbox{Minimize over $\sigma$}:&& \textbf{E}_x^\sigma\left[e^{\sum_{n=0}^\infty l(Y_n,A_n,Y_{n+1})}\right]=:\textbf{V}(x,\sigma),~x\in \textbf{X}.
\end{eqnarray}
It is also referred to as the risk-sensitive DTMDP problem.
We denote the value function of problem (\ref{GGZyExponentialProblem2}) by $\textbf{V}^\ast$. Then a strategy $\sigma^\ast$ is called optimal for problem  (\ref{GGZyExponentialProblem2}) if $\textbf{V}(x,\sigma^\ast)=\textbf{V}^\ast(x)$ for each $x\in \textbf{X}.$

\begin{condition}\label{GGZyExponentialDTMDPCon2}
\begin{itemize}
\item[(a)] The function $l(x,a,y)$ is lower semicontinuous in $a\in \textbf{A}$ for each $x,y\in \textbf{X}.$
\item[(b)] For each bounded measurable function $f$ on $\textbf{X}$ and each $x\in \textbf{X},$ $\int_{\textbf{X}}f(y)p(dy|x,a)$ is continuous in $a\in\textbf{A}.$
\item[(c)] The space $\textbf{A}$ is a compact Borel space.
\end{itemize}
\end{condition}

\begin{proposition}\label{GGZyExponentialProposition02}
Suppose Condition \ref{GGZyExponentialDTMDPCon2} is satisfied.
\begin{itemize}
\item[(a)]  The value function $\textbf{V}^\ast$ is the minimal $[1,\infty]$-valued measurable solution to \begin{eqnarray}\label{ZyExponential02}
\textbf{V}(x)=\inf_{a\in \textbf{A}}\left\{\int_{\textbf{X}}p(dy|x,a)e^{l(x,a,y)}\textbf{V}(y)\right\},~x\in \textbf{X}.
\end{eqnarray}
\item[(b)] Let $\textbf{U}$ be a $[1,\infty]$-valued lower semianalytic function on $\textbf{X}$. If
\begin{eqnarray*}
\textbf{U}(x)\ge \inf_{a\in \textbf{A}}\left\{\int_{\textbf{X}}p(dy|x,a)e^{l(x,a,y)}\textbf{U}(y)\right\},~\forall~x\in \textbf{X},
\end{eqnarray*}
then $\textbf{U}(x)\ge \textbf{V}^\ast(x)$ for each $x\in \textbf{X}.$ In particular, if the function $\textbf{U}$ satisfying the above relation is $[1,\infty)$-valued, then so is the value function $\textbf{V}^\ast.$
\item[(c)] Let $\varphi$ be a deterministic stationary strategy for the DTMDP model $\{\textbf{X},\textbf{A},p,l\}$. If
\begin{eqnarray}\label{GGZyExponential025}
\textbf{V}^\ast(x)=\int_{\textbf{X}}p(dy|x,\varphi(x))e^{l(x,\varphi(x),y)}\textbf{V}^\ast(y),~\forall~x\in \textbf{X},
\end{eqnarray}
then $\textbf{V}^\ast(x)=\textbf{V}(x,\varphi)$ for each $x\in \textbf{X}.$
\item[(d)] Let $\textbf{V}^{(0)}(x):=1$ for each $x\in \textbf{X}$, and for each $n=1,2,\dots,$
    \begin{eqnarray*}
    \textbf{V}^{(n)}(x):=\inf_{a\in A}\left\{\int_{\textbf{X}}p(dy|x,a)e^{l(x,a,y)}\textbf{V}^{(n-1)}(y)\right\},~\forall~x\in \textbf{X}.
    \end{eqnarray*}
Then $(\textbf{V}^{(n)}(x))$ increases to $\textbf{V}^\ast(x)$ for each $x\in \textbf{X}$, where $\textbf{V}^\ast$ is the value function for problem (\ref{GGZyExponentialProblem2}). Furthermore, there exists a deterministic stationary strategy $\varphi$ satisfying (\ref{GGZyExponential025}), and so in particular, there exists a deterministic stationary optimal strategy for the risk-sensitive DTMDP problem (\ref{GGZyExponentialProblem2}).
\end{itemize}
\end{proposition}

\par\noindent\textbf{Acknowledgement.} We thank the referees for their remarks, which improved the presentation of this paper. This work is partially supported by a grant from the Royal Society (IE160503).

\end{document}